\documentclass{article}

\usepackage[T1]{fontenc}
\usepackage[utf8]{inputenc}
\usepackage[british]{babel}

\usepackage{lmodern}
\usepackage[margin=1in]{geometry}

\usepackage{amsmath}
\usepackage{amsthm}
\usepackage{amsfonts}
\usepackage{hyperref}
\hypersetup{
	colorlinks = true,
	linkcolor = purple, 
	filecolor = purple, 
	citecolor = purple,
	urlcolor = purple
}
\usepackage{url} 
\usepackage{amssymb}
\usepackage{tikz-cd}
\usetikzlibrary{babel}
\allowdisplaybreaks

\usepackage[backend=biber, style=alphabetic, sorting=nyt, maxcitenames=7,
mincitenames=1, maxnames=7, abbreviate=true, block=space]{biblatex}

\newtheorem{lemma}{Lemma}[section]
\newtheorem{theorem}[lemma]{Theorem}
\newtheorem{proposition}[lemma]{Proposition}
\newtheorem{corollary}[lemma]{Corollary}
\newtheorem{conjecture}[lemma]{Conjecture}
\newtheorem{introconjecture}{Conjecture}

\newtheorem*{theorem*}{Theorem}
\newtheorem*{proposition*}{Proposition}

\theoremstyle{definition}
\newtheorem{rem}[lemma]{Remark}
\newtheorem{example}[lemma]{Example}

\newtheorem{definition}[lemma]{Definition}

\newcommand{\A}{\mathbb{A}}
\newcommand{\PP}{\mathbb{P}}

\newcommand{\Z}{\mathbb{Z}}
\newcommand{\Q}{\mathbb{Q}}
\newcommand{\R}{\mathbb{R}}
\newcommand{\C}{\mathbb{C}}

\newcommand{\G}{\mathbb{G}}

\DeclareMathOperator{\Spec}{Spec}

\DeclareMathOperator{\id}{id}

\DeclareMathOperator{\rank}{rank}
\DeclareMathOperator{\sign}{sign}
\DeclareMathOperator{\disc}{disc}

\DeclareMathOperator{\Sym}{Sym}

\DeclareMathOperator{\End}{End}
\DeclareMathOperator{\GW}{GW}

\DeclareMathOperator{\Hom}{Hom}

\DeclareMathOperator{\Tr}{Tr}

\DeclareMathOperator{\Bl}{Bl}

\DeclareMathOperator{\characteristic}{char}

\let\det\undefined
\DeclareMathOperator{\det}{det}

\newcommand{\SH}{\text{\textnormal{SH}}}

\newcommand{\Var}{\text{\textnormal{Var}}}

\newcommand{\topo}{\text{\textnormal{top}}}
\newcommand{\cat}{\text{\textnormal{cat}}}
\newcommand{\Et}{\text{\textnormal{Ét}}}
\newcommand{\EtLin}{\text{\textnormal{ÉtLin}}}
\newcommand{\Gal}{\text{\textnormal{Gal}}}
\newcommand{\Chow}{\text{\textnormal{Chow}}}

\makeatletter
\newcommand{\colim@}[2]{%
	\vtop{\m@th\ialign{##\cr
			\hfil$#1\operator@font colim$\hfil\cr
			\noalign{\nointerlineskip\kern1.5\ex@}#2\cr
			\noalign{\nointerlineskip\kern-\ex@}\cr}}%
}
\newcommand{\fcolim}{%
	\mathop{\mathpalette\colim@{\rightarrowfill@\textstyle}}\nmlimits@
}
\newcommand{\colim}{%
	\mathop{\mathpalette\colim@{}}\nmlimits@
}
\makeatother

\title{The $\A^1$-Euler characteristic of symmetric powers}
\author{Louisa F. Bröring, Jesse Pajwani, and Anna M. Viergever}
\date{}

\begin{document}
	\maketitle
	
	\begin{abstract}
		The $\mathbb{A}^1$-Euler characteristic is a refinement in algebraic geometry of the classical topological Euler characteristic, which can be constructed using motivic homotopy theory. This invariant is a quadratic form rather than an integer, which carries a lot of information, but is difficult to compute in practice. In this survey, we discuss a conjectural way for computing the $\A^1$-Euler characteristic of the symmetric powers of a variety in terms of the $\A^1$-Euler characteristic of the variety itself formulated using the theory of power structures. We discuss evidence towards the conjecture so far, techniques to approach it, and applications. 
	\end{abstract}
	
	\tableofcontents
	
	\section*{Introduction}
	\addcontentsline{toc}{section}{Introduction}
	
	The $\mathbb{A}^1$-Euler characteristic is an analogue in algebraic geometry of the classical Euler characteristic from topology. It was first introduced by Hoyois \cite{HoyoisQRGLVTF} (for smooth projective schemes), and is given by applying the categorical trace construction of Dold--Puppe \cite{dold80-categoricalEC} to the stable motivic homotopy category $\SH(k)$ of a base field $k$ which is not of characteristic $2$. Rather than being an integer, the $\mathbb{A}^1$-Euler characteristic of a variety over $k$ takes values in the  Grothendieck-Witt ring $\text{GW}(k)$ of non-degenerate quadratic forms over $k$ up to isometry. By work of Arcila-Maya, Bethea, Opie, Wickelgren, and Zakharevich \cite{Wickelgren} in characteristic zero and of Levine, Pepin Lehalleur, and Srinivas \cite{LevineLehalleurSrinivas} in general, one can view the $\mathbb{A}^1$-Euler characteristic as a motivic measure 
	\[\chi_c\colon K_0(\Var_k) \to \GW(k)\]
	where $K_0(\Var_k)$ denotes the Grothendieck ring of varieties. 
	
	The $\mathbb{A}^1$-Euler characteristic contains topological information: if $k\subset\mathbb{R}$ then the rank of $\chi_c(X)$ is equal to the compactly supported topological Euler characteristic of $X(\mathbb{C})$, while the signature of $\chi_c(X)$ is the compactly supported topological Euler characteristic of $X(\mathbb{R})$. The invariant $\chi_c$ can be used to obtain quadratic enrichments of familiar formulas like the Gau\ss -Bonnet formula by work of Déglise-Jin-Khan \cite{DegliseFCMHT} and Levine-Raksit \cite{LevineGB} and the Riemann-Hurwitz formula by work of Levine \cite{LevineA}. This invariant plays an important role in the fast-growing field of quadratically enriched enumerative geometry, which aims to obtain refinements of results in classical enumerative geometry using motivic homotopy theory. The resulting invariants, or counts, are quadratic forms that allow to recover known results over $\C$ or $\R$ and extend them over different fields. In general, $\mathbb{A}^1$-Euler characteristics are hard to compute. 
	
	In this survey, we focus on the case of symmetric powers. Write $S_n$ for the symmetric group on $n$ letters. For a CW complex $T$, write $\Sym^nT := T^n/S_n$, where $S_n$ acts on $T^n$ by permuting the factors. MacDonald \cite{MacDonald} proves that the Euler characteristic of $\Sym^nT$ satisfies
	\[
	\chi_c^{\topo}(\Sym^nT) = \binom{\chi_c^{\topo}(T) + n-1}{n}.
	\]
	
	Given a quasi-projective variety $X$ over a field $k$ not of characteristic $2$, one can ask if one can describe the $\mathbb{A}^1$-Euler characteristic of $\Sym^n(X):=X^n/S_n$ in terms of that of $X$, i.e. if one can find a quadratic refinement of the above formula.
	
	The second author and P\'al \cite{Pajwani-PalPS} formulate this problem in terms of power structures. Intuitively, a power structure on a ring $R$ is a way to make sense of the expression $f(t)^r$ for $r\in R$ and $f\in 1+tR[[t]]$ (see Definition \ref{defn:pstruc}). One can equivalently define power structures by maps $b_n\colon R\to R$ satisfying certain conditions (see Proposition \ref{proposition: power functions}). On the Grothendieck ring of varieties $\tilde{K_0}(\Var_k)$ with universal homeomorphisms (see Definition~\ref{def:universal-homeo}) inverted, there is a power structure given by taking symmetric powers, i.e. for a variety $X$ over $k$, one has $b_n([X]) = [\Sym^n(X)]$. This was first proven in characteristic zero by Gusein-Zade, Luengo and Melle-Hern\'andez \cite[Theorem~2]{Gusein-Zade-Luengo-Melle-HernandezPS}. Bejleri-McKean \cite[Corollary 5.4]{mckean2024symmetricpowersnullmotivic} show that the $\mathbb{A}^1$-Euler characteristic admits a factorization 
	\[
	\chi_c\colon K_0(\Var_k)\to \tilde{K_0}(\Var_k) \to \GW(k).
	\]
	One can then ask if there exists a power structure $b_*$ on $\GW(k)$ such that the $\A^1$-Euler characteristic respects power structures, i.e. such that $\chi_c(\Sym^n(X)) = b_n(\chi_c(X))$ for every variety $X$ over $k$ and every $n \ge 0$. In \cite{Pajwani-PalPS}, the second author and P\'al construct an explicit power structure $a_*$ on $\GW(k)$ with the property that if there exists a power structure $b_*$ on $\GW(k)$ such that $\chi_c$ respects the power structures, then $b_* = a_*$ (see Definition \ref{defn:anpstruc}). It is at present unknown whether $\chi_c$ does actually respect power structures if we equip $\GW(k)$ with the power structure $a_*$. This article discusses the following conjecture, which was also formulated by Bejleri-McKean \cite[Conjecture~1.1 and Theorem~1.2]{mckean2024symmetricpowersnullmotivic}. 
	
	\begin{introconjecture}\label{introconjecture}
		For any field $k$ of characteristic not $2$ and for any quasi-projective variety $X/k$ and $n \geq 1$
		\[
		\chi_c(\Sym^nX) = a_n(\chi_c(X)).
		\]
		In particular, $\chi_c\colon \tilde{K_0}(\Var_k) \to \GW(k)$ respects the power structures on both rings.
	\end{introconjecture}
	
	If Conjecture \ref{introconjecture} were to hold, this would give a powerful computational method to compute the $\mathbb{A}^1$-Euler characteristic of the symmetric powers of a variety. Additionally, this would be an important step towards proving a quadratically enriched analogue of the classical G\"ottsche formula, which gives a generating series for the Euler characteristic of the Hilbert schemes of points of given length of a smooth projective surface. Furthermore, little is presently known about the behavior of the $\mathbb{A}^1$-Euler characteristic under taking quotients, so understanding the case of symmetric powers may be helpful in further investigating this. 
	
	Conjecture \ref{introconjecture} has been proven in several cases. If $k$ is a field of $2$-virtual cohomological dimension zero (e.g., if $k=\R$ or $\C$), or equivalently if $\GW(k)$ is torsion free, this conjecture is largely classical, following from work of MacDonald \cite{MacDonald}, see \cite[Proposition 9.1]{mckean2024symmetricpowersnullmotivic} for a proof. Similarly if $k$ has characteristic zero with virtual cohomological dimension $1$, we may combine the above with \cite[Corollary 8.20]{Pajwani-PalYZ}.
	
	Alternatively, one can look at classes of varieties for which Conjecture \ref{introconjecture} holds. The second author and P\'al \cite{Pajwani-PalPS} show that Conjecture \ref{introconjecture} holds for all varieties of dimension zero. The second author, Rohrbach and the third author \cite{PajwaniRohrbachViergever} prove that it holds for a class of varieties which includes the classes of $[\A^n]$ and $[\mathbb{P}^n]$ and cellular varieties. The first and third author \cite{bv2024quadratic} show that Conjecture \ref{introconjecture} holds for all curves. The latter statement was independently proven by Taelman and Pepin Lehalleur. The first author \cite{broeringtorus} proves that Conjecture \ref{introconjecture} holds for all smooth projective schemes if $n=2$ and if $\characteristic(k)\neq 3$, also $n=3$. In characteristic zero, this implies that Conjecture \ref{introconjecture} holds for $n=2$ and $n=3$. 
	
	In this survey, we give an introduction to Conjecture \ref{introconjecture}, and explain more about the above-mentioned evidence towards it. In Section \ref{section: Euler characteristic introduction}, we give an introduction to the $\mathbb{A}^1$-Euler characteristic. In Section \ref{section: power structures}, we define power structures and introduce the power structure of the second author and Pál \cite{Pajwani-PalPS}. Section \ref{section: results} surveys known results towards Conjecture \ref{introconjecture}. In Section \ref{section: lambda structures}, we give an outlook on a potential strategy to prove the conjecture by lifting everything to the level of $K$-theory spectra. Finally, in Section \ref{section: Gottsche formula}, we survey work of Bejleri--McKean \cite{mckean2024symmetricpowersnullmotivic} which shows how Conjecture \ref{introconjecture} implies a quadratically enriched analogue of the classical G\"ottsche formula.
	
	\subsection*{Acknowledgements}
	\addcontentsline{toc}{subsection}{Acknowledgements}
	
	We are thankful to Marc Levine for suggesting to us to write this conference volume paper, and for numerous extremely helpful discussions about aspects of this problem. We would also like to extend our sincere thanks to the organizers of the Park City Mathematics Institute Summer Session 2024, where we gave the talk which this paper is based on, for inviting us to be there and for the wonderful research and discussion environment. Finally, we deeply thank Dori Bejleri, Stephen McKean, Ambrus P\'al, Simon Pepin Lehalleur, Herman Rohrbach and Lenny Taelman for all of their work on this problem. The first author was supported by the Deutsche Forschungsgemeinschaft
	(DFG, German Research Foundation) -- Research Training Group 2553 --
	Projektnummer 412744520 and the Faculty of Mathematics at the University of Duisburg-Essen. The second author was supported by Dan Loughran's UKRI Future Leaders' MR/V021362/1 and would like to thank the University of Bath for their support.
	
	\subsection*{Conventions}
	Throughout, let $k$ be a field of characteristic $\neq 2$. A variety over $k$ will mean a reduced separated scheme of finite type over $k$. 
	
	\section{$\A^1$-Euler characteristic}\label{section: Euler characteristic introduction}
	Let $T$ be a finite $CW$ complex. The \emph{(compactly supported) Euler characteristic} of $T$, denoted $\chi_c^{\topo}(T)$, is defined to be 
	\[
	\chi^{\topo}_c(T) := \sum_{i=0}^{\infty}(-1)^i \dim_{\Q} H^i_{sing,c}(T; \Q) \in \Z,
	\]
	where the $H^i_{sing,c}(T; \Q)$ denote the rational singular cohomology groups with compact support of $T$ with coefficients in $\Q$. This is a highly classical invariant, with its study dating back to Euler in the 18th century.
	
	The $\A^1$-Euler characteristic is an analogue of the topological Euler characteristic in algebraic geometry, which was first studied by Hoyois \cite{HoyoisQRGLVTF}, and defines a motivic measure taking values in the Grothendieck--Witt ring of quadratic forms over $k$: see \cite[Sections~2--5]{mckean2024symmetricpowersnullmotivic} for a more comprehensive introduction.
	
	\begin{definition}
		The \emph{Grothendieck-Witt ring of quadratic forms} $\GW(k)$ over $k$ is the group completion of the monoid of isometry classes of non-degenerate quadratic forms over $k$ with respect to orthogonal direct sums. The group $\GW(k)$ becomes a ring via tensor products of quadratic forms.
	\end{definition}
	
	\begin{rem}
		\label{rem:gw-gen-rel}
		As $k$ has characteristic not two, non-degenerate quadratic forms over $k$ correspond with non-degenerate symmetric bilinear forms. Via this correspondence, one can show that $\GW(k)$ is generated as an Abelian group by the one-dimensional forms $\langle a\rangle: x \mapsto ax^2$ for $a \in k^\times$, and one can find explicit relations among the $\langle a\rangle$ for $a \in k^\times$ to obtain a presentation of $\GW(k)$ (see e.g. \cite[Lemma~3.9]{MorelATF}). The form $H := \langle 1\rangle + \langle -1\rangle$ is called the \emph{hyperbolic form}.
	\end{rem}
	
	\begin{example}
		\begin{enumerate}
			\item A quadratic form over $\C$ is completely determined by its rank, so $\GW(\C) = \Z$. More generally, the rank homomorphism $\GW(k) \to \Z$ is an isomorphism whenever $k$ is quadratically closed.
			\item By Sylvester's law of inertia, a quadratic form over $\R$ is completely determined by its rank and its signature, and so $\GW(\R) \cong \Z[C_2]$ where $\Z[C_2]$ is the group ring of the cyclic group of order two. 
			\item We have $\GW(k) \cong \Z\oplus \Z/2$ as a group if $k$ is a finite field (of characteristic not two), see \cite[Chapter~2, Theorem~3.8]{ScharlauQHF}.
		\end{enumerate}
	\end{example}
	
	\begin{definition}
		The \emph{Grothendieck ring of varieties} over $k$, denoted $K_0(\Var_k)$, is the free Abelian group generated by isomorphism classes $[X]$ of varieties $X$ over $k$ modulo the scissors relation: for a variety $X$ and a closed subvariety $Z \subset X$, we have
		\[
		[X] = [Z] + [X\setminus Z] \in K_0(\Var_k).
		\]
		The group $K_0(\Var_k)$ becomes a ring via the product rule $[X]\cdot [Y] := [(X\times_k Y)_{red}]$ for $X$ and $Y$ varieties over $k$, where $(-)_{red}$ denotes taking the associated reduced closed subscheme. A ring homomorphism out of $K_0(\Var_k)$ is called a \emph{motivic measure}.
	\end{definition}
	
	\begin{rem}
		In general, $K_0(\Var_k)$ is generated by smaller classes of varieties. For example, $K_0(\Var_k)$ is generated by quasi-projective varieties, see \cite[Chapter 2, Corollary 2.6.6 (a)]{CLNS2018MotivicIntegration}. 
		
		In characteristic zero Bittner \cite[Theorem~3.1]{bittner_universal_2004} constructs a presentation of $K_0(\Var_k)$ as an Abelian group generated by the classes of connected smooth projective varieties over $k$, modulo two explicit relations: an identification of $[\emptyset]$ with $0$ and a blow-up relation. 
	\end{rem}
	
	\begin{rem}
		The ring $K_0(\Var_k)$ is not a well-behaved ring. For example, if $k$ is of characteristic zero, the class of $\A^1$ is a zero-divisor in $K_0(\Var_k)$, see Borisov \cite[Theorem 2.13]{Borisov2018ZeroDiv}. Nevertheless, it is the ring where one would like (compactly supported) Euler characteristics to live since its structure already captures many desirable properties of a (compactly supported) Euler characteristic.
	\end{rem}
	
	Following \cite{LevineLehalleurSrinivas}, we give a brief overview of the construction of the $\A^1$-Euler characteristic. This construction takes place in the stable motivic homotopy category $\SH(S)$ for a Noetherian base scheme $S$, first introduced by Voevodsky \cite{Voevodsky98HT}, and relies on the associated six functor formalism. See Hoyois \cite{HoyoisSOEMHT} for an introduction to $\SH(S)$ and its six functor formalism.
	
	Let $e$ denote the exponential characteristic of $k$. For a quasi-projective scheme $p \colon X \to \Spec k$ over $k$, we can consider the pushforward with compact support $p_!1_X \in \SH(k)$ of the unit element $1_X \in \SH(X)$. Using results of Riou when $k$ is perfect, namely \cite[Corollary B.2]{Riou} and \cite[Th\'eor\`eme 2.2]{Riou2005Duality}, and Elmanto--Khan \cite[Theorem 3.2.1]{ElmantoKhan} for a general field $k$, one can show that $p_!1_X$ is strongly dualisable in the sense of Dold and Puppe \cite[Theorem 1.3]{dold80-categoricalEC} after considering $p_!1_X$ as an element of $\SH(X)[e^{-1}]$; that is, there exists an object $(p_!1_X)^\vee \in \SH(k)[e^{-1}]$ together with morphisms $\delta_X\colon 1_k \to p_!1_X\otimes (p_!1_X)^\vee$ and $\operatorname{ev}_X\colon (p_!1_X)^\vee \otimes p_!1_X \to 1_k$ such that the compositions
	\[
	\begin{tikzcd}
		X \ar[r, "\cong", phantom] & 1_{\mathcal C}\otimes X \ar[r, "\delta_X\otimes \id"] &
		X\otimes X^\vee \otimes X \ar[r, "\id\otimes
		\operatorname{ev}_X"] &  X\otimes
		1_{\mathcal C} \ar[r, "\cong", phantom] & X
	\end{tikzcd}
	\]
	and
	\[
	\begin{tikzcd}
		X^\vee \ar[r, "\cong", phantom] & X^\vee\otimes 1_{\mathcal C} \ar[r, "\id\otimes\delta_X"] &
		X^\vee \otimes X\otimes X^\vee \ar[r, "\operatorname{ev}_X\otimes\id"] &
		1_{\mathcal C}\otimes X^\vee \ar[r, "\cong", phantom] & X^\vee
	\end{tikzcd}
	\]
	are the identity. The triple $((p_!1_X)^\vee, \delta_X, \operatorname{ev}_X)$ is unique up to unique isomorphism. For a strongly dualisable object $p_!1_X$, define its categorical Euler characteristic $\chi^{\text{cat}}(p_!1_X) \in \End_{\SH(k)[e^{-1}]}(1_k)$ as the composition
	\[
	\begin{tikzcd}
		1_k \ar[r, "\delta_X"] &
		p_!1_X\otimes (p_!1_X)^\vee \ar[r, "\cong", phantom] & (p_!1_X)^\vee \otimes p_!1_X \ar[r, "\operatorname{ev}_X"] &
		1_k
	\end{tikzcd}
	\]
	where the middle isomorphism is the switching isomorphism.
	
	One can now extend a computation of Morel \cite[Corollary~1.24]{MorelATF}, or for $k$ non-perfect of Bachmann-Hoyois \cite[Theorem~10.12]{BachmannHoyoisNMHT}, to obtain an isomorphism $\End_{\SH(k)}(1_k)[e^{-1}] \cong \GW(k)[e^{-1}]$.
	
	\begin{definition}
		The \emph{(compactly supported) $\A^1$-Euler characteristic of a quasi-projective $k$-scheme $X$} with structure map $p\colon X \to \Spec k$, denoted $\chi_c(X)$ is defined as the image of $\chi^\cat(p_!1_X)$ in $\GW(k)[e^{-1}]$.
	\end{definition}
	
	One can now compute that the rank of $\chi_c(X)$ agrees with the compactly supported étale Euler characteristic of $X$ and hence is an integer. Thus, $\chi_c(X)$ lies in the subring $\GW(k) \subset \GW(k)[e^{-1}]$ and we can consider $\chi_c(X)$ to be an element of $\GW(k)$. One can now show \cite[46]{LevineLehalleurSrinivas} that $\chi_c$ is compatible with the scissors relation and products, and thus furnishes a ring homomorphism $\chi_c: K_0(\Var_k) \to \GW(k)$.
	\begin{definition}
		The motivic measure 
		\[\chi_c\colon K_0(\Var_k) \to \GW(k)\]
		constructed above is called the \emph{(compactly supported) $\A^1$-Euler characteristic}. 
	\end{definition}
	
	If $X$ is a smooth projective $k$-scheme, there is a more hands-on description of its $\A^1$-Euler characteristic, which is a consequence of Levine--Raksit's motivic Gauß-Bonnet theorem.
	
	\begin{theorem}[Motivic Gauß-Bonnet, Levine--Raksit, {\cite[Corollary 8.7]{LevineGB}}]
		\label{thm:motivicGB}
		Let $X$ be a connected, smooth, projective variety over $k$. Then there exists an integer $m \in \Z$ controlled by the dimensions of the Hodge cohomology groups $H^q(X,\Omega^p_X)$ for $p,q \in \Z$ such that
		\begin{itemize}
			\item If $\dim X$ is odd, then $\chi_c(X) = m\cdot H \in \GW(k)$.
			\item If $\dim X = 2n$ is even, then $\chi_c(X) = \beta + m\cdot H \in \GW(k)$ where $\beta$ is the quadratic form corresponding with the symmetric bilinear form
			\[
			\begin{tikzcd}
				H^n(X,\Omega^n_X) \times H^n(X, \Omega^n_X) \ar[r, "\cup"] & H^{2n}(X,\Omega^{2n}_X) \ar[r, "\Tr_{X/k}"] & k
			\end{tikzcd}
			\]
			where $\cup$ denotes the cup product in Hodge cohomology and $\Tr_{X/k}$ denotes the trace from Serre duality.
		\end{itemize}
	\end{theorem}
	
	\begin{rem}
		We can use Theorem \ref{thm:motivicGB} to give an equivalent definition of the $\A^1$-Euler characteristic if $k$ has characteristic zero. By Bittner's presentation \cite[Theorem~3.1]{bittner_universal_2004}, we only need to define $\chi_c$ for connected smooth projective varieties and verify that this definition satisfies the scissors relation for blow-ups. Therefore, we can define $\chi_c(X)$ for connected smooth projective $k$-varieties $X$ using Theorem \ref{thm:motivicGB} and then use \cite[Product formula on page 2184 and Proposition 2.4 (5)]{LevineA}
		to verify that $\chi_c(X)$ satisfies the product and blow up relations, so extends to a motivic measure. Alternatively one can explicitly verify the product relation and the scissors relation for blow-ups via an explicit computation using Hodge cohomology, see \cite[Proposition~5.31 and Proposition~5.35]{broeringtorus}. This construction of the compactly supported $\A^1$-Euler characteristic in characteristic zero predates the categorical definition and is due to Arcila-Maya, Bethea, Opie, Wickelgren, and Zakharevich \cite[Theorem 2.13]{Wickelgren}.
	\end{rem}
	
	\begin{example}
		\label{chi-c-examples}
		\begin{enumerate}
			\item Let $L/k$ be a finite, separable field extension. Then $\chi_c(\Spec L) = [\Tr_{L/k}] \in \GW(k)$, where $\Tr_{L/k}$ denotes the class in $\GW(k)$ of the field trace bilinear form on $L$: $(x,y) \mapsto \Tr_{L/k}(xy)$. This was first shown by Hoyois, see \cite[Theorem 1.9]{HoyoisQRGLVTF}. 
			\item $\chi_c(\PP^n_k) = \sum_{i=0}^n\langle (-1)^i\rangle =
			\begin{cases}
				\langle 1\rangle + \frac{n}{2} \cdot H & \text{for }n \text{ even},\\
				\frac{n+1}{2}\cdot H & \text{for }n \text{ odd}.
			\end{cases}$ 
			\item $\chi_c(\A^n_k) = \langle (-1)^n\rangle$.
			\item $\chi_c(\G_m) = \langle -1\rangle - \langle 1\rangle$.
		\end{enumerate}
		We note that (2) was first shown by Hoyois, see \cite[Example 1.7]{HoyoisQRGLVTF}. Examples (3) and (4) above are straightforward computations from the $\mathbb{P}^n_k$ case, using that $\chi_c$ is a ring homomorphism.
	\end{example}
	
	By work of Levine \cite[Remark 2.3 (1)]{LevineA}, the $\A^1$-Euler characteristic relates to the topological Euler characteristic of real and complex points.
	\begin{proposition}[Levine, {\cite[Remark 2.3 (1)]{LevineA}}]\label{prop: relation chi to topological Euler characteristics}
		Let $X$ be a smooth, projective, $k$-variety. Let $\chi_c^\topo$ denote the compactly supported topological Euler characteristic. Then
		\begin{enumerate}
			\item If $k \subset \C$, then $\rank \chi_c(X) = \chi_c^\topo(X(\C))$.
			\item If $k \subset \R$, then the signature of $\chi_c(X)$ with respect to the embedding $k \subset \R$ agrees with $\chi_c^\topo(X(\R))$.
		\end{enumerate}
		Moreover, by viewing $\chi_c$ as a motivic measure, we can remove the assumption that $X$ is smooth and projective (see \cite[Theorem 2.14 and Theorem 2.20]{Pajwani-PalYZ} for details).
	\end{proposition}
	Another important observation that is important in the study of the $\A^1$-Euler characteristic of symmetric powers is that the $\A^1$-Euler characteristic is insensitive to universal homeomorphisms
	
	\begin{definition}\label{def:universal-homeo}
		A morphism of schemes $f\colon X \to Y$ is a \emph{universal homeomorphism} if $X\times_YY' \to Y'$ is a homeomorphism of underlying topological spaces for all morphisms of schemes $Y' \to Y$.
	\end{definition}
	
	\begin{definition}
		Consider the ideal
		\[
		I = ([X]-[Y]\mid \text{there is a universal homeomorphism } f\colon X\to Y) \subset K_0(\Var_k)
		\]
		generated by differences of universally homeomorphic $k$-varieties. We define the \emph{Grothendieck ring of $k$-varieties up to universal homeomorphism} as $\tilde K_0(\Var_k) = K_0(\Var_k)/I$. Note that when $\mathrm{char}(k)=0$, $I$ is the zero ideal and $K_0(\Var_k) = \tilde K_0(\Var_k)$.
	\end{definition}
	
	\begin{proposition}[Bejleri--McKean, {\cite[Corollary 5.4]{mckean2024symmetricpowersnullmotivic}}]
		The $\A^1$-Euler characteristic factors over $\tilde K_0(\Var_k)$, that is if $X$ and $Y$ are universally homeomorphic, then $\chi_c(X) = \chi_c(Y)$.
	\end{proposition}
	
	The $\A^1$-Euler characteristic behaves nicely with respect to Zariski- and Nisnevich-locally trivial fibre bundles, see for example \cite[Section 5]{mckean2024symmetricpowersnullmotivic}. Unlike the topological Euler characteristic, the $\A^1$-Euler characteristic does not behave well with respect to étale locally trivial fibre bundles. This complicates the study of the $\A^1$-Euler characteristic of quotients: indeed even if a variety $X$ has a nice, free action by a nice finite group $G$, we do not in general have a simple description of $\chi_c(X/G)$ in terms of $\chi_c(X)$ and $\chi_c(G)$: see for example \cite[Proposition 5.11]{mckean2024symmetricpowersnullmotivic} for details.
	
	\section{Power structures}\label{section: power structures}
	
	\subsection{Definition} 
	Let $R$ be a commutative ring. A power structure on $R$ is a structure which mimics taking powers of a power series over $\Z$. More precisely, one can make the following definition due to Gusein-Zade, Luengo and Melle-Hernández \cite{Gusein-Zade-Luengo-Melle-HernandezPS,Gusein-Zade-Luengo-Melle-HernandezHS}.
	
	\begin{definition}\label{defn:pstruc}
		Let $R$ be a commutative ring. A \emph{power structure} on $R$ is a map
		\[
		(1+tR[[t]])\times R\to (1+tR[[t]]), (f(t),r)\mapsto f(t)^r
		\]
		such that:
		\begin{itemize}
			\item $f(t)^0=1$ and $f(t)^1=f(t)$
			\item $(f(t)g(t))^r = f(t)^rg(t)^r$
			\item $f(t)^{r+s} = f(t)^rf(t)^s$
			\item $f(t)^{rs}=(f(t)^r)^s$
			\item $(1+t)^m= 1+mt+O(t^2)$
			\item If $g(t)=f(t)^m$ then $f(t^n)^m=g(t^n)$ for any $n\in\mathbb{Z}_{\geq 0}$
			\item (Finite generation property:) For any $N > 0$, there exists $M > 0$ such that if
			$f (t) \in 1 + t \cdot R[[t]]$ and $m\in R$, we may determine $f(t)^m\mod t^N$ solely from $m$ and $f (t)\mod t^M$. 
		\end{itemize}
	\end{definition}
	
	Given a commutative ring $R$ equipped with a power structure $b_*$, for $r\in R$, define $b_i(r)$ to be such that 
	\[
	(1-t)^{-r} = 1 + \sum_{i=0}^\infty b_i(r)t^i.
	\]
	By the last three properties above, every $f(t)\in 1+R[[t]]$ can be written as $f(t) = \prod_{i=0}^\infty (1-t^i)^{-c_i}$ for $c_i\in R$, which in turn means that the $b_i$ uniquely determine the power structure. The following proposition utilises this to give an equivalent definition of power structures which is often more useful in practice. 
	
	\begin{proposition}[Gusein-Zade--Luengo--Melle-Hernández, {\cite[Proposition 1]{Gusein-Zade-Luengo-Melle-HernandezHS}}]\label{proposition: power functions}
		On a ring $R$, functions $b_i\colon R\to R$ for $i\in\mathbb{Z}_{\geq 0}$ such that:
		\begin{itemize}
			\item $b_0(r) = 1$ for all $r$
			\item $b_1(r)=r$ for all $r$
			\item $b_i(0) = 0$ and $b_i(1) =1$ for all $i\geq 1$
			\item $b_n(r+s) = \sum_{i=0}^nb_i(r)b_{n-i}(s)$
		\end{itemize}
		uniquely extend to a power structure on $R$.
	\end{proposition}
	
	\begin{example}
		On $\mathbb{Z}$, there is a power structure given by defining $f(t)^n = f(t) \cdots f(t)$ for $n \ge 0$ and $f(t)^n = (\frac{1}{f(t)})^{-n}$ for $n < 0$. Equivalently, we can describe this power structure as
		\[(1-t)^{-n} = \sum_{i=0}^\infty \binom{n+i-1}{n-1} t^i .\]
	\end{example}
	
	\begin{example}\label{example: power structure on K0}
		On the Grothendieck group of varieties with universal homeomorphisms inverted, $\tilde{K_0}(\Var_k)$, there is a power structure given by 
		\[
		(1-t)^{-[X]} = 1 + \sum_{i=0}^{\infty} [\Sym^n(X)] t^n
		\]
		whenever $X/k$ is a quasiprojective variety. One can check this using \cite[Lemma 7.29]{MustataZF}. In characteristic zero, this was first proven by Gusein-Zade, Luengo and Melle-Hern\'andez \cite[Theorem~2]{Gusein-Zade-Luengo-Melle-HernandezPS}, using the
		standard stratification of the symmetric power. If $\operatorname{char}(k)>0$, the standard stratification of the $n^{\text{th}}$ symmetric power of a variety $X$ only becomes isomorphic to $\Sym^n(X)$ after inverting universal homeomorphisms, which is why we need to consider the power structure on $\tilde{K_0}(\Var_k)$. 
	\end{example}
	
	\begin{definition}
		Equip $\tilde{K_0}(\Var_k)$ with the power structure from the above example. By a slight abuse of notation, for any $c \in \tilde{K_0}(\Var_k)$ write $\Sym^n(c)$ for the coefficient of $t^n$ in the power series $(1-t)^{-c}$. We see that if $X/k$ is a quasiprojective variety then $\Sym^n([X]) = [\Sym^n(X)]$.
	\end{definition}

	\begin{definition}\label{definition: map respects power structures}
		Suppose that $R$ and $R'$ carry power structures. Let $\varphi\colon R\to R'$ be a ring morphism, and
		\[
		\tilde{\varphi}\colon R[[t]]\to R'[[t]], \sum_{i=0}^\infty c_it^i\mapsto \sum_{i=0}^\infty \varphi(c_i)t^i
		\]
		be the induced map. We say that $\varphi$ \emph{respects the power structure} if $\tilde{\varphi}(f(t)^a) = \tilde{\varphi}(f(t))^{\varphi(a)}$ for all $f(t)\in 1 + tR[[t]]$ and $a\in R$. In terms of the $b_i$ functions, if we write $b_i^R$ and $b_i^{R'}$ for the respective functions defining the power structures, it is easy to check that $\varphi$ \emph{respects the power structure} if 
		\[
		\varphi(b_i^R(r)) = b_i^{R'}(\varphi(r))
		\] 
		for all $i\in\mathbb{Z}_{\geq 0}$ and $r\in R$.
	\end{definition}
	
	\subsection{Euler characteristics of symmetric powers}
	By Example \ref{example: power structure on K0}, the ring $\tilde{K_0}(\Var_k)$ carries a power structure given by symmetric powers. One can wonder if there is a power structure $b_*$ on $\GW(k)$ such that the $\A^1$-Euler characteristic $\chi_c\colon \tilde{K_0}(\Var_k)\to \GW(k)$ respects the power structures in the sense of Definition \ref{definition: map respects power structures}, i.e. such that $\chi_c(\Sym^n([X])) = b_n(\chi_c(X))$ for all varieties $X$. Working towards this, we recall some classical results on topological Euler characteristics. 
	
	Let $T$ be a CW complex. Define $T^{(n)} :=T^n / S_n$, where the symmetric group $S_n$ acts on $T^n$ by permuting the factors. The main theorem of MacDonald \cite[Page 568]{MacDonald} implies the following.
	
	\begin{theorem}[MacDonald, {\cite[Page 568]{MacDonald}}]\label{thm:MacDonald}
		The Euler characteristic $\chi_c^\topo(T^{(n)})$ is given by ${ \chi_c^\topo(T) + n - 1 \choose n}$.
	\end{theorem}
	\begin{corollary}\label{cor: conjecture for algebraically closed field}
		Suppose $k$ is a separably closed field, so that $\GW(k) \cong \Z$ and let $X/k$ be a variety. Then $\chi_c(X^{(n)}) = { \chi_c(X) + n - 1 \choose n}$. 
	\end{corollary}
	\begin{proof}
		If $k$ is of characteristic $0$, we may use a standard field of definition argument, see \cite[Theorem 2.14]{Pajwani-PalYZ}, to reduce to the case where $k=\mathbb{C}$. Here we may identify $\chi_c(X) = \chi_c^{\topo}(X(\mathbb{C}))$ by Proposition \ref{prop: relation chi to topological Euler characteristics}, which therefore reduces to Theorem \ref{thm:MacDonald}. Furthermore, Bejleri-McKean \cite[Proposition 4.1]{mckean2024symmetricpowersnullmotivic} use a comparison theorem between the $\ell$-adic Euler characteristic and the rank of $\chi_c$ to extend the result to positive characteristic.
	\end{proof}
	
	\begin{corollary}\label{cor:RealMacDonald}
		Let $k=\mathbb{R}$. Then $\sign( \chi_c(\mathrm{Sym}^n(X)))$ is the coefficient of $t^n$ in the power series
		\[
		(1-t)^{-\chi_c^{\topo}(X(\mathbb{R}))} (1-t^2)^{- (\chi_c^{\topo}(X(\mathbb{C}))-\chi_c^{\topo}(X(\mathbb{R})))/2} .
		\]
	\end{corollary}
	\begin{proof}
		This essentially follows instantly by Corollary~\ref{cor:RealMacDonald} and a computation of $\mathrm{Sym}^n(X)(\mathbb{R})$ due to Kharlamov and Rasdeaconu (see \cite[\S3.2]{KharlamovRasdeaconu}).
	\end{proof} 
	
	\subsection{Power structures on $GW(k)$}
	
	The rich structure of quadratic forms as well as the wealth of field arithmetic at our disposal allows us to equip $\GW(k)$ with much additional structure that goes beyond the ring structure. Here, we survey some known results about power structures on $\GW(k)$. Vector spaces over $k$ come equipped with a natural notion of a symmetric power: namely, $S^nV := V^{\otimes n}/S_n$, where the symmetric group acts by permuting the tensors. Following McGarraghy in \cite{McGarraghySP}, it is natural to ask whether there are natural notions of ``symmetric powers'' of a quadratic form. This was answered in part by McGarraghy~\cite[Definition~3.1]{McGarraghySP}.
	\begin{definition}
		Let $V$ be a vector space over $k$ and let $\varphi\colon V \times V \to k$ be a non-degenerate symmetric bilinear form on $k$. Then the \emph{factorial symmetric power} of $\varphi$ is a symmetric bilinear form $\mathbf{S}^n\varphi$ on $S^nV$, defined on elementary tensors by
		
		\[
		\mathbf{S}^n\varphi( (x_1 \odot \ldots \odot x_n), (y_1 \odot \ldots \odot y_n)) = \operatorname{per}( \varphi(x_i,y_j))_{1\leq i,j \leq n},
		\]
		where $\operatorname{per}$ denotes the permanent of an $n \times n$ matrix (see \cite[Definition~3.1]{McGarraghySP}) and $x_1 \odot \ldots \odot x_n$ is the image of $x_1 \otimes \ldots \otimes x_n$ under the quotient map $V^{\otimes n} \to S^nV$.
	\end{definition}
	Now, \cite[Proposition 3.5]{McGarraghySP} ensures that if $[\varphi] = \sum_{i=1}^m \langle a_i \rangle \in \GW(k)$, then 
	\[
	[\mathbf{S}^n\varphi] = \sum_{\substack{1 \leq i_1 \leq i_l \leq m:\\n_{i_1} + \ldots + n_{i_l} =n}} \prod_{j=1}^l \langle n_{i_j}! a_{i_j}^{n_{i_j}} \rangle.
	\]
	In particular, if $\operatorname{char}(k) = p$ and $n>p$, then $\mathbf{S}^n\varphi$ is degenerate. To remedy this, McGarraghy introduced the \emph{non-factorial symmetric power of $\varphi$} to be the power structure $S^\bullet$ on $\GW(k)$ defined by
	\[
	S^n(\langle \alpha \rangle) = \langle \alpha^n \rangle.
	\]
	This has nicer algebraic properties on the ring $\GW(k)$, but given a bilinear form $\varphi\colon  V \times V \to k$, there is no canonical bilinear form $S^n\varphi\colon S^nV \times S^nV \to k$ such that $[S^n\varphi] = S^n([\varphi]) \in \GW(k)$; constructions of this quadratic form involve choosing a basis.
	
	As in the paper by the second author and P\'al (\cite[Corollary~3.26]{Pajwani-PalPS}), we define another power structure on $\GW(k)$ using functions as in Propostion \ref{proposition: power functions}. 
	
	\begin{definition}\label{defn:anpstruc}
		Define the power structure $a_*$ on $\GW(k)$ by 
		\begin{equation}\label{equation: power structure from Pajwani-Pal}
			a_n(\langle \alpha \rangle) := \langle \alpha^n\rangle + \frac{n(n-1)}{2}(\langle 2 \rangle + \langle \alpha \rangle - \langle 1\rangle - \langle 2\alpha\rangle ).
		\end{equation}
	\end{definition}
	A priori, this does not appear like a natural power structure, however it is justified by the fact that $a_*$ can be characterized as the only power structure on $\GW(k)$ such that if $X/k$ is a variety of dimension $0$, we have that $\chi_c(\Sym^nX) = a_n(\chi_c(X))$ (see \cite[Corollary~3.27]{Pajwani-PalPS}).

	The term $\langle 2 \rangle + \langle \alpha \rangle - \langle 1\rangle - \langle 2\alpha\rangle$ is $2$-torsion in $\GW(k)$, see \cite[Corollary~3.18]{Pajwani-PalPS}. Indeed, it can be rewritten as $\langle \langle 2, \alpha \rangle \rangle - 2\mathbb{H}$, where $\langle \langle 2,\alpha \rangle \rangle$ denotes the Pfister form. Writing $[\alpha]$ to mean the class represented by $\alpha$ in $H^1(k, \Z/2\Z) \cong k^\times/k^{\times2}$, we obtain the following result, either using the isomorphism from the Milnor conjectures, or using an elementary argument as in \cite{Pajwani-PalPS}.
	\begin{lemma}[{\cite[Lemma 3.29]{Pajwani-PalPS}}]
		The term $\langle 2 \rangle + \langle \alpha \rangle - \langle 1\rangle - \langle 2\alpha\rangle$ vanishes if and only if $[2] \cup [\alpha] = 0$ as an element of $H^2(k, \Z/2\Z)$.
	\end{lemma}
	
	\begin{rem}\label{rem: non-vanishing of ta}
		As noted in \cite[Corollary 3.30]{Pajwani-PalPS}, this means that the above power structure agrees with McGarraghy's non-factorial symmetric powers whenever the map $[2] \cup - : H^1(k, \Z/2\Z) \to H^2(k, \Z/2\Z)$ is $0$. This is true whenever $k$ is a field such that $2$ is a square, or if $\mathrm{vcd}_2(k)\leq 1$ (see Definition \ref{defn: vcd}).
	\end{rem}
	\begin{rem}
		The constructions above from \cite{McGarraghySP} do not just construct a power structure on $\GW(k)$. These power structures have slightly more data: namely, if $\varphi: V \times V \to k$ is a bilinear form, then we have a (not necessarily canonical) bilinear form $S^n \varphi: S^n V \times S^n V \to k$ such that $[S^n \varphi] = S^n([\varphi])$.
		
		This is no longer true when we look at the power structure $a_\ast$ from \cite{Pajwani-PalPS}. For example, let $\alpha \in k^\times$ be such that $[2] \cup [\alpha] \neq 0 \in H^2(k, \Z/2\Z)$. Then $a_2(\langle \alpha \rangle) = \langle 1 \rangle + \langle 2 \rangle + \langle \alpha \rangle - \langle 1 \rangle - \langle 2\alpha \rangle$ is not an effective element of $\GW(k)$. That is, there is no bilinear form $\varphi: V \times V \to k$ such that $[\varphi] = a_2(\langle \alpha \rangle)$. The power structure we have constructed therefore does not give rise to a notion of ``symmetric power of bilinear forms''; it only gives us a notion of a power structure on $\GW(k)$.
	\end{rem}
	
	\section{Results}\label{section: results}
	
	Endowing $\tilde{K_0}(\Var_k)$ with the power structure from Example \ref{example: power structure on K0} and $\GW(k)$ with the power structure $a_*$ from Equation \eqref{equation: power structure from Pajwani-Pal}, we can now view the $\A^1$-Euler characteristic $\chi_c\colon \tilde{K_0}(\Var_k)\to \GW(k)$ as a map between rings that carry a power structure. We have now introduced everything to formulate the following precise version of Conjecture \ref{introconjecture}, which was also formulated by Bejleri-McKean \cite[Conjecture~1.1 and Theorem~1.2]{mckean2024symmetricpowersnullmotivic}.
	
	\begin{conjecture}\label{conjecture: main conjecture}
		For any field $k$ of characteristic not $2$ and for any quasi-projective variety $X/k$ and $n \geq 1$
		\[
		\chi_c(\Sym^nX) = a_n(\chi_c(X)).
		\]
		In particular, $\chi_c$ respects the power structures in the sense of Definition \ref{definition: map respects power structures} if we endow $\tilde{K_0}(\Var_k)$ with the power structure from Example \ref{example: power structure on K0} and $\GW(k)$ with the power structure $a_*$.
	\end{conjecture} 
	
	If this conjecture were true, the equality $\chi_c(\Sym^n(X)) = a_n(\chi_c(X))$ would let us readily compute $\chi_c(\Sym^n(X))$ for any $n\ge 0$ from $\chi_c(X)$. 
	The power structure $a_*$ is the only power structure that one can hope to use for such a computational purpose, due to the following result. 
	\begin{theorem}[{\cite[Corollary~3.27]{Pajwani-PalPS}}]\label{thm:uniquenessofan}
		Endow $\tilde K_0(\Var_k)$ with the power structure from Example \ref{example: power structure on K0}.
		Suppose that there exists a power structure $b_*$ on $\GW(k)$ such that $\chi_c$ respects power structures after endowing $\GW(k)$ with the power structure $b_*$, then $b_* = a_*$.
	\end{theorem}
	
	\begin{proof}[Proof sketch]
		Suppose that a power structure $b_*$ exists with $\chi_c(\Sym^n(c)) = b_n(\chi_c(c))$ for all $c \in \tilde{K_0}(\Var_k)$. Since elements of the form $\langle \alpha \rangle$ generate $\GW(k)$ as an Abelian group as in Remark \ref{rem:gw-gen-rel}, it is enough to show that $b_n(\langle \alpha \rangle)=a_n(\langle \alpha \rangle)$ for any $\alpha \in k^\times/(k^{\times})^2$. 
		
		Taking $c = \Spec(k(\sqrt{\alpha}))$ yields $\chi_c(c) = \langle 2 \rangle + \langle 2\alpha \rangle$. Computing $\Sym^n(c)$ explicitly yields $\Sym^n(c) = \frac{n+1}{2}[c]$ for $n$ odd and $\frac{n}{2}[c] + [\Spec(k)]$ for $n$ even. By our hypothesis that $b_*(\chi_c(c))=\chi_c(\Sym^n(c))$, we have the following equality whenever $n$ is odd:
		\[
		\frac{n+1}{2}(\langle 2 \rangle + \langle 2\alpha \rangle) = \chi_c(\Sym^n(\Spec(k(\sqrt{\alpha}))) = b_n(\langle 2 \rangle + \langle 2\alpha \rangle) = \sum_{i=0}^n b_i (\langle 2 \rangle) b_{n-i}(\langle 2\alpha\rangle),
		\]
		and an analogous expression when $n$ is even. Proceeding by induction on $n$ and specialising to $2\alpha=1$ implies $b_n(\langle2\rangle)=a_n(\langle 2 \rangle)$, which we may in turn use to see that $b_n(\langle \alpha \rangle)=a_n(\langle \alpha \rangle)$.
	\end{proof}
	There are two possible approaches to proving Conjecture \ref{conjecture: main conjecture}. The first is to show that the conjecture holds for certain classes of fields. In light of this, it is helpful to have the following definition.
	\begin{definition}\label{defn: vcd}
		Write $\operatorname{cd}_2(k)$ for the $2$-cohomological dimension of $k$ (see \cite[\S3.1]{SerreGaloisCohomology}). We say $k$ has \emph{$2$-virtual cohomological dimension $n$} and write $\operatorname{vcd}_2(k)=n$ if there exists a finite extension $L/k$ with $\operatorname{cd}_2(L) = n$. 
	\end{definition}
	Note that $\operatorname{vcd}_2(k)=n$ is equivalent to $\operatorname{cd}_2(k(\sqrt{-1}))=n$ by Artin--Schreier theory (i.e., combining \cite[Satz 3]{ArtinSchreier} with \cite[Théor{\`e}me]{SerreCD}). The following result is essentially a corollary of Theorem \ref{thm:MacDonald}.
	\begin{corollary}\label{corollary: 2vcd 0}
		Suppose that $k$ is a field such that $\operatorname{vcd}_2(k)=0$. Then Conjecture \ref{conjecture: main conjecture} holds over $k$.
	\end{corollary}
	\begin{proof}
		The conjecture holds if we can show that for any variety $X/k$, we have $\chi_c(\Sym^n(X))=a_n(\chi_c(X))$.
		
		Now $k$ has $2$-virtual cohomological dimension $0$. If $k$ does not have $2$-cohomological dimension $0$, then by Artin--Schreier theory (see e.g. \cite[Proposition 9.8.1]{Scheiderer}), $k$ admits a real closure $k^{rc}$. That is, $k^{rc}$ is an algebraic extension of $k$ which is not algebraically closed but $k^{rc}(\sqrt{-1})$ is. Moreover, this forces $\operatorname{char}(k)=0$. 
		
		By \cite[Chapter 2, \S2, Proposition 2]{SerreGaloisCohomology}, there are no non-trivial quadratic extensions of $k(\sqrt{-1})$. Therefore since $\GW(k)$ is generated by classes $\langle \alpha \rangle$ for $\alpha \in k^\times/k^{\times2}$, we see that $\GW(k)$ is either generated solely by $\langle 1 \rangle$ if $k$ has $2$-cohomological dimension $0$, or $\langle 1 \rangle$ and $\langle -1 \rangle$ else. 
		
		If $k = k(\sqrt{-1})$, then the base change map $\GW(k) \to \GW(\overline{k})$ is an isomorphism so the statement follows by base changing to the algebraic closure and applying Corollary \ref{cor: conjecture for algebraically closed field}, where we use the equality between $\rank(\chi_c(X))$ and $\chi_c^{\topo}(X(\mathbb{C}))$ (or the $\ell$-adic Euler characteristic in positive characteristic).
		
		If $k \neq k(\sqrt(-1))$, then the base change map $\GW(k) \to \GW(k^{rc})$ is an isomorphism. In the case where $k^{rc}=\mathbb{R}$, Corollary \ref{cor:RealMacDonald} gives the result. For $k^{rc}$ a general real closed field, we may argue as in \cite[Theorem 2.20]{Pajwani-PalYZ} and use techniques of Pál \cite[Definition 2.5]{PalRealClosedSimplicial} to reduce the argument to the real case. Note that this result uses the equality between $\sign(\chi_c(X))$ and the real closed  Euler characteristic $\chi_c^{\topo}(X(k))$.
	\end{proof}
	\begin{rem}
		The Milnor conjectures allow us to see that the condition $\operatorname{vcd}_2(k) =0$ is equivalent to $\GW(k)$ being torsion free. In this language, an alternative proof of the above statement appears as \cite[Corollary~9.2]{mckean2024symmetricpowersnullmotivic}. We choose to use the language of cohomological dimension since we also have the following theorem, which essentially appears as \cite[Lemma~8.20]{Pajwani-PalYZ}.
	\end{rem}
	\begin{theorem}\label{Theorem: vcd 1}
		Suppose that $k$ has characteristic $0$ and $\operatorname{vcd}_2(k) = 1$. Then Conjecture \ref{conjecture: main conjecture} holds over $k$.
	\end{theorem}
	One example of a field $k$ as in the above theorem would be if $k$ is a function field in one variable over $\mathbb{R}$. 
	\begin{proof}
		Since $\operatorname{cd}_2(k(\sqrt{-1}))\leq1$, the Milnor conjectures (proven by Voevodsky in \cite{VoevodskyMilnorConjecture}) imply that an element of $\GW(k(\sqrt{-1}))$ is completely determined by its rank and discriminant, and thus, an element of $\GW(k)$ is determined by its rank, signature and discriminant. Therefore, we need to show that if $X/k$ is a variety over $k$, and $\varphi$ is either $\rank, \sign$ or $\disc$, then $\varphi(\chi_c(\Sym^n(X))) = \varphi(a_n(\chi_c(X)))$. For $\varphi=\rank$ or $\sign$, this directly follows by Corollary \ref{corollary: 2vcd 0}. We therefore only need to show that $\disc(\chi_c(\Sym^n(X))) = \disc(a_n(\chi_c(X)))$ for every variety $X/k$. Therefore fix a variety $X/k$, and without loss of generality, using the same methods as \cite[proof of Theorem 2.14]{Pajwani-PalYZ}, we may assume that $k$ is finitely generated over $\Q$.
		
		Let $q \in \GW(k)$. An elementary computation yields
		\[
		\disc(a_n(q)) = \disc(q)^{ {n + \rank(q) - 1 \choose n}},
		\]
		where we first show the result for elements of the form $\langle \alpha \rangle$, and then use an induction to obtain the results for general elements of $\GW(k)$. The result therefore holds if
		\[
		\disc(\chi_c(\Sym^n(X))) = \disc(\chi_c(X))^{{ n + \rank(\chi_c(X)) - 1 \choose n}}.
		\]
		
		For any variety $X/k$, define $w(X) := \frac{ \rank(\chi_c(X)) \cdot \dim(X)}{2}$. Let $\Gal_k$ denote the absolute Galois group of $k$. Note that $\disc(\chi_c(\Sym^n(X))) \in k^\times/(k^{\times})^2$, which we identify with $\Hom(\Gal_k, \Z/2\Z)$ by Kummer theory since $k$ is not of characteristic~$2$. Fix a prime $\ell$, so that embedding $\Z/2\Z = \{\pm 1\} \hookrightarrow \Q_\ell^\times$ allows us to view $\disc(\chi_c(\Sym^n(X)))$ as a $1$-dimensional $\ell$-adic representation of $\Gal_k$. Similarly, we obtain a $1$-dimensional $\ell$-adic representation of $\Gal_k$, $\det_\ell(X)$, defined by
		\[
		\det\nolimits_\ell(X) := \prod_{i=0}^{2 \dim(X)} H^i_c(X_{\overline{k}}, \Q_\ell)^{\otimes (-1)^i}.
		\]
		Since $k$ has characteristic $0$, we may apply \cite[Theorem 2.29]{Pajwani-PalYZ}, which uses Saito's Theorem \cite[Theorem 2]{Saito}, to obtain that for any variety $X/k$ of dimension $d$
		\[
		\disc(\chi_c(X)) \cdot (-1)^{ w(X)} = \Q_\ell( -w(X)) \cdot \det_\ell(X) \in \Hom(\Gal_k, \Q_\ell^\times).
		\]
		As in \cite[Remark 2.30]{Pajwani-PalYZ}, since $k$ is finitely generated over $\Q$, the cyclotomic character has infinite order so $w(X)$ is the unique integer such that  $\Q_\ell( -w(X)) \cdot \det_\ell(X)$ is torsion. Therefore to show the required equality between $\disc(\chi_c(\Sym^n(X)))$ and $\disc(a_n(\chi_c(X)))$, it is enough to show the corresponding equality of the $\mathrm{det}_\ell$ representations, i.e.
		\[
		\det\nolimits_\ell( \Sym^n(X)) = \det\nolimits_\ell(X)^{ {n + \rank\chi_c(X) - 1 \choose n}}.
		\]
		We therefore compute $\det\nolimits_\ell(\Sym^n(X))$. Note that $\det_\ell(X)$ is the determinant representation $\det(H^*_c( \Sym^n(X)_{\overline{k}}, \Q_\ell))$, where here $\det$ denotes the graded determinant as in \cite[Notation 2.41]{Pajwani-PalYZ}.  Consider the $\Q_\ell$-graded-commutative algebra $H^*_c( \Sym^n(X)_{\overline{k}}, \Q_\ell)$. By the symmetric K\"unneth formula, \cite[Expos\'e XVII, Th\'eor\`eme 5.5.21]{SGA4}, we obtain
				\[
		H^*_c( \Sym^n(X)_{\overline{k}}, \Q_\ell)  = \Sym^n( H^*_c( X_{\overline{k}}, \Q_\ell)),
		\]
		where $\Sym^n(H^*_c( X_{\overline{k}}, \Q_\ell))$ is the $n^{\text{th}}$ symmetric power of this graded-commutative algebra, see  \cite[Expos\'e XVII, Th\'eor\`eme 5.5.21]{SGA4}. A slight alteration of \cite[Lemma 8.20]{Pajwani-PalYZ} to the graded setting gives us
		\[
		\det( \Sym^n( H^*_c( X_{\overline{k}}, \Q_\ell))) = \det( H^*_c(X_{\overline{k}}, \Q_\ell))^{{ n + \rank\chi_c(X) - 1 \choose n}},
		\]
		which yields the result. 
	\end{proof}

	\begin{rem}\label{remark: higher stiefel whitney classes}
		These results all use comparisons between invariants extracted from $\chi_c(X)$ and other classical cohomological invariants: the case where $k=\mathbb{C}$ uses the equality $\rank(\chi_c(X)) = \chi_c^{\topo}(X(\mathbb{C}))$, the real case uses that $\sign(\chi_c(X)) = \chi_c^{topo}(X(\mathbb{R}))$, and the cohomological dimension $1$ case uses that $\disc(\chi_c(X))$ can be expressed in terms of the compactly supported $\ell$-adic cohomology $H^*_{\acute{e}t, c}(X_{\overline{k}}, \Q_\ell)$ (see for example, \cite[Theorem 2]{Saito}, \cite[Theorem 2.24]{Pajwani-PalYZ}), even when $X$ is not smooth and projective. The difficulty with this approach is that we do not understand the higher cohomological invariants coming from $\chi_c(X)$. Saito offers a conjectural formula \cite[Conjecture 2.3]{SaitoHW} comparing the Hasse--Witt invariant of $\chi_c(X)$ to an invariant coming from $\ell$-adic cohomology. However, this is still a conjecture. To the knowledge of the authors, there are no current conjectures relating higher Stiefel--Whitney classes of $\chi_c(X)$ to other invariants.
		
		It should also be noted that whenever $\mathrm{vcd}_2(k)\leq 1$ the power structure above agrees with McGarraghy's non-factorial symmetric power structure: see Remark \ref{rem: non-vanishing of ta}. It is an open question whether Conjecture \ref{conjecture: main conjecture} holds in any cases where this is not true.
	\end{rem}
	
	Bejleri and McKean \cite[Theorem~7.10]{mckean2024symmetricpowersnullmotivic} show that Conjecture \ref{conjecture: main conjecture} for all varieties is equivalent to showing Conjecture \ref{conjecture: main conjecture} for varieties whose $\A^1$-Euler characteristic is zero. Examples of these are elliptic curves, for which the conjecture holds by \cite[Theorem~27]{bv2024quadratic}. More generally, the $\A^1$-Euler characteristic of an Abelian variety is known to be zero by \cite[Theorem 5.29]{Pajwani-PalYZ}, but it is, at present, unknown if this also holds for their symmetric powers.
	
	Aside from the above results, which seek to prove the full conjecture for a given field $k$, there are also many results which instead focus on enlarging the class of varieties $X$ for which $\chi_c(\Sym^n(X)) = a_n(\chi_c(X))$ over an arbitrary field $k$. In light of this, it is helpful to have the following definition.
	\begin{definition}
		Let $\tilde{K_0}(\Sym_k)$ denote the subset of $\tilde{K_0}(\Var_k)$ given by
		\[
		\tilde{K_0}(\Sym_k) := \{ c \in \tilde{K_0}(\Var_k): \chi_c(\Sym^n(c)) = a_n(\chi_c(c)) \text{ for all $n$}\}.
		\]
	\end{definition}
	In this language, Conjecture \ref{conjecture: main conjecture} is equivalent to showing $\tilde{K_0}(\Sym_k)=\tilde{K_0}(\Var_k)$.

	The papers \cite{Pajwani-PalPS}, \cite{PajwaniRohrbachViergever} and \cite{bv2024quadratic} all approach this conjecture by showing that for certain varieties $X$ you have $[X] \in \tilde{K_0}(\Sym_k)$, with the goal being to eventually show that either $\tilde{K_0}(\Sym_k) = \tilde{K_0}(\Var_k)$, or explicitly find a counterexample. There are two major approaches to showing a variety is contained in $\tilde{K_0}(\Sym_k)$, which we outline in the rest of this section.

	\subsection{Algebraic properties of $\tilde{K_0}(\Sym_k)$}
	The first approach is to prove algebraic properties about the subset $\tilde{K_0}(\Sym_k)$. It is clear that $0 \in \tilde{K_0}(\Sym_k)$ and $[\Spec(k)] \in \tilde{K_0}(\Sym_k)$. Moreover, a quick computation yields the following.
	\begin{lemma}[{\cite[Lemma 2.9]{Pajwani-PalPS}}]\label{Lemma: K0Sym is an Abelian subgroup}
		The set $\tilde{K_0}(\Sym_k)$ is an Abelian subgroup of $\tilde{K_0}(\Var_k)$. That is, if $r,s \in \tilde{K_0}(\Sym_k)$, then $-r, r+s \in \tilde{K_0}(\Sym_k)$.
	\end{lemma}
	
	\begin{definition}
		Let $K_0(\Et_k)$ denote the subring of $\tilde{K_0}(\Var_k)$ generated by varieties of dimension $0$. Similarly, let $K_0(\EtLin_k)$ denote the subring of $\tilde{K_0}(\Var_k)$ generated by $[\A^1]$ and $K_0(\Et_k)$.
	\end{definition}
	The main result of \cite{Pajwani-PalPS} yields the following structure on $\tilde{K_0}(\Sym_k)$. 
	\begin{theorem}[{\cite[Corollary 4.30]{Pajwani-PalPS}}]\label{Theorem: K0Sym is a K0Et submodule}
		The subset $\tilde{K_0}(\Sym_k)$ is a $K_0(\Et_k)$-submodule of $\tilde{K_0}(\Var_k)$. In particular, $K_0(\Et_k)\subset \tilde{K_0}(\Sym_k)$.
	\end{theorem}
	\begin{proof}[Sketch of proof]
		The proof here proceeds in four steps. We first claim that $[\Spec(k(\sqrt{\alpha}))] \in \tilde{K_0}(\Sym_k)$ for any $\alpha \in k^\times$. This is essentially immediate by Theorem \ref{thm:uniquenessofan}. 
		
		The second step is to show that $[\Spec(k(\sqrt{\alpha}, \sqrt{\beta}))] \in \tilde{K_0}(\Sym_k)$ for any $\alpha, \beta \in k^\times$. This goes via a direct computation of $\Sym^n(\Spec(k(\sqrt{\alpha}, \sqrt{\beta})))$.
		
		The third step is to use the above result to show that if $X/k$ is a variety with $[X] \in \tilde{K_0}(\Sym_k)$, then $[X] \cdot [\Spec(k(\sqrt{\alpha}))] \in \tilde{K_0}(\Sym_k)$. This essentially proceeds by a Galois descent argument as in \cite[Theorem 3.8]{Pajwani-PalYZ}, by noting that $[X \amalg X] \in \tilde{K_0}(\Sym_k)$ and that $\Sym^n(X \amalg X)$ is a $k(\sqrt{\alpha})/k$ form of $\Sym^n(X \times_k \Spec(k(\sqrt{\alpha})))$, allowing us to define descent data on $\chi_c(\Sym^n(X \amalg X))$ to recover $\chi_c(\Sym^n(X \times_k \Spec(k(\sqrt{\alpha})))$. This allows us to reduce to the case where $X$ is a quadratic extension of $k$, and then the result follows by the biquadratic case from step 2. Repeatedly applying this result, we see that if $c \in \tilde{K_0}(\Sym_k)$, then if $A/k$ is any multiquadratic étale algebra, $c \cdot [A] \in \tilde{K_0}(\Sym_k)$.

		The fourth and final step is to use a powerful result of Garibaldi--Merkurjev--Serre \cite[Theorem 29.1]{GMS}. Fix $c \in \tilde{K_0}(\Sym_k)$. Then the assignment $A \mapsto \chi_c( \Sym^n( c \cdot [A] )) - a_n (\chi_c(c \cdot [A])) \pmod{H}$ gives us a cohomological invariant valued in $\operatorname{W}(k)$, where $\operatorname{W}(k) = \GW(k)/(H)$ is the \emph{Witt ring}. The previous step shows that this invariant is $0$ whenever $A/k$ is a multiquadratic étale algebra. 
		
		In \cite[Theorem 29.1]{GMS}, it is shown that if a Witt ring valued cohomological invariant vanishes for all multiquadratic étale algebras, then it vanishes for all finite étale algebras. In particular, this implies that $\chi_c( \Sym^n( c \cdot [A] )) = a_n (\chi_c(c \cdot [A])) \pmod{H}$ for any $c \in \tilde{K_0}(\Sym_k)$. to get the result in $\GW(k)$, we then only need to prove that these sides have the same rank, which follows by a basic computation. In particular, if $c \in \tilde{K_0}(\Sym_k)$, then $c \cdot [A] \in \tilde{K_0}(\Sym_k)$ for any finite étale algebra $A$, which yields the result.
	\end{proof}
	\begin{corollary}
		Conjecture \ref{conjecture: main conjecture} holds whenever $X/k$ is a variety of dimension $0$.
	\end{corollary}
	\begin{proof}
		Using that $1 = [\Spec(k)] \in \tilde{K_0}(\Sym_k)$ and that $[X] \in K_0(\Et_k)$, the above theorem implies that $[X] \cdot 1 = [X] \in \tilde{K_0}(\Sym_k)$.
	\end{proof}
	
	Combining the above result with G\"ottsche's Lemma for symmetric powers of $\A^n$ (see \cite[Lemma~4.4]{GottscheLemma} or \cite[Chapter~7, Proposition~1.1.11]{CLNS2018MotivicIntegration}), the second author, Rohrbach and the third author prove the following.
	
	\begin{theorem}[{\cite[Theorem~4.10 and Corollary~4.12]{PajwaniRohrbachViergever}}]
		The subset $\tilde{K_0}(\Sym_k)$ is a $K_0(\EtLin_k)$-submodule of $\tilde{K_0}(\Var_k)$. In particular $K_0(\EtLin_k)\subset \tilde{K_0}(\Sym_k)$.
	\end{theorem}
	
	\begin{proof}[Sketch of proof.]
		G\"ottsche's Lemma for symmetric powers \cite[Lemma 4.4]{GottscheLemma} or \cite[Proposition 1.1.11 in Chapter 7]{CLNS2018MotivicIntegration} states that for $X$ a variety over $k$, we have an equality $[\Sym^n(X\times \A^m)] = [\Sym^n X][\A^{mn}]$ in $\tilde{K_0}(\Var_k)$. Taking $X = [\Spec(k)]$ shows that $[\A^n]\in \tilde{K_0}(\Sym_k)$. Using Theorem~\ref{Theorem: K0Sym is a K0Et submodule}, we find that Conjecture~\ref{conjecture: main conjecture} holds for varieties in  $K_0(\EtLin_k)$. To see that $\tilde{K_0}(\Sym_k)$ is a $K_0(\EtLin_k)$-submodule of $\tilde{K_0}(\Var_k)$, we use that for $[X]\in \tilde{K_0}(\Sym_k)$, we have that $\chi_c([\Sym^n(X\times \A^m)]) =  a_n(X\times \A^m)$ by G\"ottsche's Lemma and a direct computation. As every element of $K_0(\EtLin_k)$ can be written in the form $\sum_{i=1}^k m_i[\A^{n_i}][\Spec (L_i)]$ where $L_i/k$ is a finite \'etale algebra, we can now deduce the result from Lemma \ref{Lemma: K0Sym is an Abelian subgroup} and Theorem \ref{Theorem: K0Sym is a K0Et submodule}. 
	\end{proof}
	
	\begin{example}
		The power of the above theorem is that the classes of many natural varieties lie in the ring $K_0(\EtLin_k)$.  For example, $[\A^n]$ and $[\mathbb{P}^n]$ are both in $K_0(\EtLin_k)$ for all $n\geq 0$. More generally, cellular varieties in the sense of \cite[Page 2189]{LevineA} are in $K_0(\EtLin_k)$.
		One can also show that the ``axis cross'' $\{xy = 0\}\subset \A^2$ defines an element of $K_0(\EtLin_k)$. Similarly, the class of the cuspidal curve $\{xy^2 = z^3\}\subset\mathbb{P}^3$ is in $K_0(\EtLin_k)$. Using the orbit-cone correspondence for split toric varieties, see for example \cite[Proposition in Section~3.1]{Fulton1993IntroToric}, the first author \cite[Proposition~3.14]{broeringtorus} shows that the classes of split toric varieties lie in $K_0(\EtLin_k)$.
	\end{example}
	
	\subsection{Direct computations}
	
	The class of varieties in $K_0(\EtLin_k)$ is still far away from every variety: indeed, every connected variety in $K_0(\EtLin_k)$ is geometrically stably rational by \cite[Theorem~6.5]{PajwaniRohrbachViergever}, so curves of genus $g\geq 1$ are outside of $\tilde{K_0}(\Var_k)$. One can nevertheless prove that $\chi_c$ respects power structures for curves, which was done by the first and third author \cite[Theorem~27]{bv2024quadratic} and independently by Taelman and Pepin Lehalleur as part of work in progress on the surface case. 
	
	\begin{theorem}
		\label{thm:curves}
		Conjecture \ref{conjecture: main conjecture} is true for all curves. 
	\end{theorem}
	
	\begin{proof}[Sketch of proof.] The proof relies on the fact that if $C$ is a smooth projective curve over $k$, then the symmetric powers of $C$ are smooth. Using the motivic Gauß-Bonnet Theorem, one can then show that for $C$ a smooth projective curve of genus $g$ and $n\in\mathbb{Z}_{\geq 0}$,  if $n = 2m$ is even, we have
		\[
		\chi_c(\Sym^{n}(C)) = \sum_{i=0}^m\binom{g}{i}\langle -1\rangle^i + \frac{1}{2}\left(\binom{2g-2}{n} - \sum_{i=0}^m\binom{g}{i}\right)H\in \GW(k)
		\]
		and if $n$ is odd, we have
		\[
		\chi_c(\Sym^{n}(C)) = -\frac 12 \binom{2g-2}{n} H\in \GW(k).
		\]
		The result for all curves is deduced from this using a normalization argument.
	\end{proof} 
	
	The first author \cite[Theorem~8.3 and Theorem~8.6]{broeringtorus} extends this method to compute $\chi_c(\Sym^2X)$ and $\chi_c(\Sym^3X)$ for a smooth, projective variety $X$ over $k$. 
	\begin{theorem}
		\label{thm:compatibility-2-3}
		Conjecture \ref{conjecture: main conjecture} holds for smooth projective schemes and $n = 2$ or if $\characteristic k \ne 3$, also $n = 3$.
	\end{theorem}
	
	Using Bittner's presentation, we can now deduce, see \cite[Theorem~8.9]{broeringtorus}:
	
	\begin{theorem}
		If $\characteristic k = 0$, then Conjecture \ref{conjecture: main conjecture} holds for $n =2$ and $n =3$.
	\end{theorem}
	
	\begin{proof}[Proof sketch]
		The idea of the proof of Theorem \ref{thm:compatibility-2-3} is similar to the proof of Theorem \ref{thm:curves}: one first computes $\chi_c(\Sym^nX)$ explicitly and then compares the computed quadratic form to $a_n(\chi_c(X))$. Unlike in Theorem \ref{thm:curves}, for general $X$, we note that $\mathrm{Sym}^n(X)$ is no longer smooth for $n \geq 2$, so to compute $\chi_c(\Sym^nX)$ we first need to resolve the singularities of $\Sym^nX$ and then compute its Hodge cohomology. For $n =2$, we can blow-up the diagonal $\Delta \subset X^2$ to obtain the scheme $\Bl_\Delta X^2$ with an $S_2$-action. The quotient by $S_2$ is again smooth and projective and an abstract blow-up of $\Sym^2X$ along $X$.
		
		One can now relate the $\A^1$-Euler characteristics of $\Bl_\Delta X^2$ and $(\Bl_\Delta X^2)/S_2$ by computing invariants in Hodge cohomology using the motivic Gauß-Bonnet theorem, and then in turn relate $\chi_c(\Bl_\Delta X^2)$ and $\chi_c((\Bl_\Delta X^2)/S_2)$ to $\chi_c(X^2)$ and $\chi_c(\Sym^2X)$ respectively using blow-up formulae. In order to keep track of the changes in the Hodge cohomology and the $S_2$ action, the first author extends the theory of an equivariant quadratic Euler characteristic first developed by the second author and Pál \cite[Section~3 in Version 2 on arXiv]{Pajwani-PalYZ}, see \cite[Section 5]{broeringtorus} for details. With this, one obtains an explicit formula for $\chi_c(\Sym^2X)$, which one can then compare to $a_2(\chi_c(X))$.
		
		The proof in the case $n = 3$ is essentially the same, but requires more than one blow-up to obtain a smooth $S_3$-quotient, which complicates the combinatorics involved in the computation. 
	\end{proof} 
	It is non-trivial to extend these computations to $n \ge 4$: for such symmetric powers, we may no longer apply the same methods as for $n \leq 3$, since we can no longer obtain a smooth quotient by blowing up in an invariant smooth subscheme.
	
	\section{Power structures and \texorpdfstring{$\lambda$}{λ}-structures}\label{section: lambda structures}
	In this section, we outline a potential approach to Conjecture \ref{conjecture: main conjecture} using higher homotopy theory, as well as highlighting some problems this approach faces.
	
	Zakharevich \cite{Zakharevich} and Campbell \cite{Campbell} give constructions of a spectrum $K(\Var_k)$ in the stable homotopy category $\SH$ with  $\pi_0(K(\Var_k)) = K_0(\Var_k)$, see \cite[Proposition 1.3]{Campbell}. In \cite[Theorem~1.1]{Nanavaty}, Nanavaty shows that if $k$ admits resolution of singularities, the $\A^1$-Euler characteristic lifts to a map of spectra in $\SH$:
	\[
	\tilde{\chi}\colon K(\Var_k)\to \End_{\SH(k)}(1_{\SH(k)})
	\]
	which recovers $\chi_c$ when taking $\pi_0$. One can wonder if one could approach Conjecture \ref{conjecture: main conjecture} by lifting to the spectrum level. To formulate this more precisely, it is helpful to have the following definition.
	
	\begin{definition}
		\label{def:lambda-structure}
		Let $R$ be a commutative ring. A \emph{pre-lambda structure} on $R$ is a collection of maps $\lambda_n\colon R\to R$ for $n\in\mathbb{Z}_{\geq 0}$ such that 
		\begin{enumerate}
			\item $\lambda_0(r) = 1$  and $\lambda_1(r) = r$ for all $r\in R$. 
			\item $\lambda_n(r+s) = \sum_{i=0}^n\lambda_i(r)\lambda_{n-i}(s)$ for all $r,s\in R$. 
			\item $\lambda_n(0) = 0$ for all $n\geq 1$
			\item $\lambda_n(1) = 0$ for all $n\geq 2$. 
		\end{enumerate}
		A pre-$\lambda$-structure is called a \emph{$\lambda$-structure} if in addition, the following axioms are satisfied:
		\begin{enumerate}\setcounter{enumi}{4}
			\item We have that 
			\[\lambda_n(xy) = P_n(\lambda_1(x),\dots, \lambda_n(x), \lambda_1(y),\dots, \lambda_n(y)).\]
			Here, $P_n(s_1^x,\dots, s_n^x, s_1^y,\dots, s_n^y)$ is the coefficient of $t^n$ in the expansion of $\prod_{i,j=1}^n (1+tx_iy_j)$ in terms of the elementary symmetric polynomials $s_i^x$ and $s_j^y$ in degree $n$ of the $x_i$ and $y_j$. 
			\item $\lambda_n(\lambda_m(x)) = P_{n,m}(\lambda_1(x),\dots, \lambda_{mn}(x))$. Here, $P_{n,m}(s_1^x,\dots, s_{nm}^x)$ is the coefficient of $t^n$ in the expansion of $\prod_{ 1\leq i_1 < \cdots < i_m\leq mn} (1+tx_{i_1},\dots, x_{i_m})$ in terms of the elementary symmetric polynomials $s_i^x$ of the $x_i$ in degree $nm$. 
		\end{enumerate}
	\end{definition}
	
	\begin{rem}\label{remark axioms 6 and 7}
		Following Larsen and Lunts \cite[Section 4]{Larsen-Lunts}, one way to think about axiom (5) is the following. Let $R$ be a commutative ring endowed with a pre-$\lambda$ structure. Consider the map 
		\[
		\varphi: R\to 1+tR[[t]], r\mapsto \sum_{n=0}^\infty \lambda_n(r)t^n.
		\]
		Note that $1+tR[[t]]$ is an Abelian group, with the usual product of power series as its group operation. By the relations (1) and (2), $\varphi$ is a morphism of groups. We can endow $1+tR[[t]]$ with a ring structure which is characterised by the identity
		\[
		\left(\prod_{i=0}^n (1+r_it) \right) \odot \left(\prod_{j=0}^m (1+s_jt) \right) = \prod_{i=0}^n \prod_{j=0}^m(1+r_is_jt). 
		\]
		The relation (5) is now equivalent to $\varphi$ being a morphism of rings, i.e. to $\varphi(rs) = \varphi(r)\odot\varphi(s)$. 
	\end{rem}
	
	From a power structure, one can obtain a pre-$\lambda$-structure and the other way around. 
	
	\begin{proposition}
		\label{prop:opposite-lambda-structure}
		Let $R$ be a commutative ring endowed with a power structure. Then setting $\lambda_n(r)$ to be the coefficient of $t^n$ in $(1+t)^r$ defines a pre-$\lambda$ structure on $R$. 
	\end{proposition} 
	
	\begin{definition}
		In the situation of Proposition \ref{prop:opposite-lambda-structure}, we call $\lambda_n$ the \emph{opposite pre-$\lambda$ structure} associated to the power structure on $R$. 
	\end{definition} 
	\begin{rem}
		When we define a power structure using functions $b_n$, these functions are defined so that $(1-t)^{-r} = \sum_{n=0}^\infty b_n(r)t^n$. Therefore, from a pre-$\lambda$ structure, we may return to a power structure: if $R$ is a ring with a pre-$\lambda$ structure, then we can construct functions $b_n: R \to R$ such that
		\[
		\sum_{n=0}^{\infty} b_n(r) t^n = \left(\sum_{n=0}^{\infty} \lambda_n(r) (-t)^n\right)^{-1},
		\]
		and it is relatively straightforward to check that these functions extend to a power structure.
	\end{rem}
	
	One of the reasons $\lambda$-structures on rings are studied is that they naturally arise from Adams operations on $K$-theory spectra. The relationship between power structures on $K_0(\Var_k)$ and Adams operations was explored by Gorsky in \cite{Gorsky} (note: Gorsky refers to a pre-$\lambda$ structure simply as a ``$\lambda$-structure'', and refers to $\lambda$-structures as ``special $\lambda$-structures'').  On page 310 of loc. cit., Gorsky shows that from a $\lambda$-structure on a ring $R$, we may construct ring homomorphisms $\Psi_n: R \to R$ called Adams operations (named due to their resemblence to Adams operations on $K$-theory spectra). One could therefore wonder if one can attempt to prove Conjecture \ref{conjecture: main conjecture} on the spectrum level via the following. First, show that there exist natural Adams operations on $K(\Var_k)$ and $\End_{\SH(k)}(1_k)$ such that the pre-$\lambda$-structure on the connected components $K_0(\Var_k)$ and $\GW(k)$ are the opposite pre-$\lambda$-structures to the power structures given by symmetric powers and the power structure $a_*$, respectively. Once we have these Adams operations, we only need to show that $\chi\colon  K(\Var_k) \to \End_{\SH(k)}(1_k)$ preserves these Adams operations. This fails at the first hurdle.
	
	\begin{proposition}
		The opposite pre-$\lambda$-structure $\lambda^{\Sym}$ of the power structure on $K_0(\Var_k)$ given by symmetric powers does not extend to a $\lambda$-structure if $k=\mathbb{C}$.
	\end{proposition}
	
	\begin{proof}
		Let $C$ be a smooth, projective curve with strictly positive genus. By classical results of Kapranov, the zeta function $\zeta_C(t) = \sum_{n=0}^\infty [\Sym^n(C)]t^n$ is a rational function. On the other hand, by \cite[Theorem~7.6]{Larsen-Lunts}, the zeta function  $\zeta_{C\times C}(t)$ is irrational in the sense of \cite[Section 2]{Larsen-Lunts}. This means that there exists a motivic measure to a field for which the image series is not a rational function. Following Remark \ref{remark axioms 6 and 7}, we consider the map  
		\[
		\varphi\colon K_0(\Var_k)\to 1 + tK_0(\Var_k)[[t]], [X]\mapsto \zeta_X(t^{-1}). 
		\]
		We have that $\lambda^{\Sym}$ satisfies axiom (6) of Definition \ref{def:lambda-structure} if and only if $\varphi$ is a morphism of rings. 
		
		By \cite[Lemma 2.5]{Larsen-Lunts}, we have that if a function is rational, then the inverse is also rational. This implies that $\zeta_{C\times C}(-t)^{-1}$ is again irrational. Therefore, there exists a motivic measure $\mu\colon K_0(\Var_k)\to K$ to a field $K$ such that $\mu(\zeta_{C\times C}(-t)^{-1})$ is not a quotient of polynomials in $K[[t]]$. 
		
		By \cite[Proposition 2.4 and Lemma 2.5]{Larsen-Lunts}, $\mu(\zeta_C(-t)^{-1})$ is a rational function. As under $\odot$, two rational functions are taken to a rational function, we see that $\mu(\zeta_C(-t)^{-1}\odot \zeta_C(-t)^{-1})$ is also rational. Therefore, $\mu(\zeta_C(-t)^{-1}\odot \zeta_C(-t)^{-1})\neq \mu(\zeta_{C\times C}(-t)^{-1})$, which implies that $\varphi([C\times C])\neq \varphi([C])\odot \varphi([C])$ and so axiom (6) of Definition \ref{def:lambda-structure} fails. 
	\end{proof}  
	\begin{corollary}
		Let $k$ be a field which embeds into $\mathbb{C}$. Then $\lambda^{\Sym}$ does not extend to a $\lambda$-structure. 
	\end{corollary}
	\begin{proof}
		Let $C$ be a smooth, projective curve over $k$ of strictly positive genus. The base change map $\psi\colon  K_0(\Var_k)\to K_0(\Var_{\mathbb{C}})$ takes the Kapranov zeta function over $k$ to the Kapranov zeta function over $\mathbb{C}$. This implies that $\zeta_C(-t)^{-1}\odot \zeta_C(-t)^{-1}$ and $\zeta_{C\times C}(-t)^{-1}$ have different images under $\psi$, and so $\zeta_C(-t)^{-1}\odot \zeta_C(-t)^{-1}\neq \zeta_{C\times C}(-t)^{-1}$ in $K_0(\Var_k)$, which means that again, axiom (6) of Definition \ref{def:lambda-structure} fails. 
	\end{proof}
	\begin{corollary}
		If $k \subseteq \mathbb{C}$, then there are no Adams operations on $K(\Var_k)$ which, on the level of $K_0$, correspond to the pre-$\lambda$-structure $\lambda^{\Sym}$.
	\end{corollary}
	\begin{proof}
		Adams operations on a ring spectrum will induce a $\lambda$-structure on the level of $\pi_0$, and $\lambda^{\Sym}$ is not a $\lambda$-structure by the above corollary.
	\end{proof}
	The above suggests that approaching the conjecture using Adams operations on $K(\Var_k)$ is impossible. Nevertheless, while the desired operation may not exist on the $K$-theory spectrum of varieties, a theorem of Heinloth shows the following.
	
	\begin{theorem}[Heinloth, {\cite[Lemma 4.1]{Heinloth}}]
		Let $K_0(\Chow_k)$ denote the Grothendieck ring of Chow motives. Then there is a $\lambda$-ring structure on $K_0(\Chow_k)$ such that the associated opposite power structure is the one induced by symmetric powers of Chow motives.
	\end{theorem}
	
	This is particularly relevant in the case of studying the $\mathbb{A}^1$-Euler characteristic $\chi_c$ due to the following.
	
	\begin{definition}
		Let $K_0(\Var_k)^{\Chow}$ denote the smallest subring of $K_0(\Chow_k)$ such that
		\begin{enumerate}
			\item The image of the natural map $K_0(\Var_k) \to K_0(\Chow_k)$ given by $X \mapsto (X, \Delta_X)$ lands in $K_0(\Var_k)^{\Chow}$.
			\item For any $M \in K_0(\Chow_k)$, the Tate twist $M(1)$ lies in $K_0(\Chow_k)$. 
		\end{enumerate}
	\end{definition}
	\begin{lemma}[{\cite[Definition 8.16]{Pajwani-PalYZ}}]
		Let $k$ be a field of characteristic $0$. Let $I := \ker(\rank\colon \GW(k) \to \Z)$, let $T$ denote the torsion subgroup of $\GW(k)$, and let $J := I^2 \cap T$. Then there is a unique ring homomorphism $\chi^{Chow}: K_0(\Var_k)^{Chow} \to \GW(k)/J$ such that the following diagram is commutative
		\begin{center}
			\begin{tikzcd}
				K_0(\Var_k) \ar[r] \ar[d, "\chi_c"] & K_0(\Var_k)^{\Chow} \ar[d, "\chi^{\Chow}"] \\
				\GW(k) \ar[r, twoheadrightarrow] & \GW(k)/J.
			\end{tikzcd}
		\end{center}
	\end{lemma}
	The above therefore sketches a roadmap to proving Conjecture \ref{conjecture: main conjecture} using $K$-theory spectra: if one could show that the map $\chi^{\Chow}$ actually extends to a map $K_0(\Chow_k) \to \GW(k)$, it may be possible to lift this to a map of spectra $K(\Chow_k) \to \End_{\SH(k)}(1_k)$. If it were possible to construct Adams operations on $K(\Chow_k)$ that realise the symmetric power $\lambda$-structure on the level of $K_0$, and to then show that these Adams operations are compatible with this morphism, then this would result in Conjecture \ref{conjecture: main conjecture}. 
	
	One difficulty of constructing the map $K_0(\Chow_k) \to \GW(k)$ is that the construction of $\chi^{Chow}$ in \cite{Pajwani-PalYZ} utilises the comparison between Stiefel--Whitney classes of $\chi_c(X)$ and constructions related to Weil cohomology theories on $X$, which are not known for higher Stiefel--Whitney classes in general (see Remark~\ref{remark: higher stiefel whitney classes}).
	
	\section{G{\"o}ttsche's Formula}\label{section: Gottsche formula}
	Bejleri and McKean \cite[\S8]{mckean2024symmetricpowersnullmotivic} apply the above work on power structures to prove a G{\"o}ttsche formula. Let $X/k$ be a smooth, projective surface. Write $X^{[n]}$ for the Hilbert scheme of points of length $n$ on $X$ and, for $x$ a closed point of $X$, write $X^{[n]}_x$ for the Hilbert scheme of points of length $n$ supported only on  $x$. We adopt the convention that $X^{[0]} = \Spec(k)$. G{\"o}ttsche \cite[Theorem~0.1]{Gottsche} proves the following beautiful formula.
	\begin{theorem}[G{\"o}ttsche, {\cite[Theorem~0.1]{Gottsche}}]\label{thm:GottscheFormula}
		
		When $k=\mathbb{C}$, there is an equality of power series in $\Z[[t]]$
		\[
		\sum_{n=0}^\infty \chi_c^{\topo}(X^{[n]}(\mathbb{C})) t^n = \prod_{n=1}^{\infty}( 1-t^n)^{-\chi_c^{\topo}(X(\mathbb{C}))}.
		\]
	\end{theorem}
	The above formula is one of the key steps in proving the Yau--Zaslow formula for counting rational curves on $K3$ surfaces, see Beauville \cite{Beauville} for more details. While not the original proof, Gusein-Zade, Luengo and Melle-Hernández provide a proof of Theorem \ref{thm:GottscheFormula} in the spirit of this paper.
	\begin{theorem}[Gusein-Zade--Lueno--Melle-Hernández, {\cite[Theorem 1]{Gusein-Zade-Luengo-Melle-HernandezHS}}]
		The power structure on $\tilde{K_0}(\Var_k)$ given by symmetric powers induces the equality
		\[
		\sum_{n \geq 0} [X^{[n]}] t^n = \sum_{n\geq 0} ( (\mathbb{A}^2)^{[n]}_0 t^n)^{[X]}.
		\]
	\end{theorem}
	Bejleri and McKean \cite[Corollary 8.4]{mckean2024symmetricpowersnullmotivic} compute $\chi_c$ of the right hand side as follows:
	\[
	\sum_{n \geq 0}  \chi_c([(\mathbb{A}^2)^{[n]}_0]) t^n = \prod_{n \geq 1} (1- \langle -1 \rangle^{n-1} t^n )^{-1} \in \GW(k)[[t]],
	\]
	Combining this computation with the above formula yields
	\begin{corollary}[Bejleri--McKean, {\cite[Theorem 8.6]{mckean2024symmetricpowersnullmotivic}}]
		Let $k$ be a field such that $\chi_c: \tilde{K_0}(\Var_k) \to \GW(k)$ preserves the power structures on both sides. Then
		\[
		\sum_{n \geq 0} \chi_c([X^{[n]}]) t^n = \prod_{n \geq 1} (1- \langle -1 \rangle^{n-1} t^n )^{-\chi_c(X)}.
		\]
	\end{corollary}
	Specialising to $k=\mathbb{C}$ yields G{\"o}ttsche's original formula, Theorem \ref{thm:GottscheFormula}. Proving Conjecture \ref{conjecture: main conjecture} would yield a G{\"o}ttsche formula over any field $k$ which is not of characteristic $2$. Moreover, work as in \cite{Pajwani-PalYZ} derives a refinement of the Yau--Zaslow formula using a quotient of $\chi_c$, obtaining a power series valued in a quotient of $\GW(k)$ which counts the rational curves on the $K3$ surface with additional arithmetic information. The above G{\"o}ttsche formula shows that proving Conjecture \ref{conjecture: main conjecture} is one step to refining the main theorem from \cite{Pajwani-PalYZ} to obtain a full arithmetic Yau--Zaslow formula over $\GW(k)$.
	
	\printbibliography[heading=bibintoc]

	%	\newpage 
	\noindent Louisa F. Br\"oring \\
	Universit\"at Duisburg-Essen \\
	Fakult\"at f\"ur Mathematik, Thea-Leymann-Str. 9, 45127 Essen, Germany \\
	E-Mail: \href{mailto:louisa.broering@uni-due.de}{louisa.broering@uni-due.de}\\
	
	\noindent Jesse Pajwani \\
	University of Bristol \\
	School of Mathematics, School of Mathematics, University of Bristol, Bristol, BS8 1TW, UK, and the Heilbronn Institute for Mathematical Research, Bristol, UK\\
	Email: \href{mailto:jesse.pajwani@bristol.ac.uk}{jesse.pajwani@bristol.ac.uk}\\
	
	\noindent Anna M. Viergever \\
	Leibniz Universit\"at Hannover\\
	Institute of Algebraic Geometry, Welfengarten 1, 30167 Hannover, Germany \\
	E-Mail: \href{mailto:viergever@math.uni-hannover.de}{viergever@math.uni-hannover.de}\\
	
	\noindent Keywords: motivic homotopy theory, refined enumerative geometry, other fields\\ 
	Mathematics Subject Classification: 14G27, 14N10, 14F42
\end{document}